\documentclass[12pt]{amsart}

\usepackage[latin1]{inputenc}
\usepackage[T1]{fontenc}      
\usepackage{amssymb,url}

\theoremstyle{plain}
\newtheorem{theorem}{Theorem}[section]
\newtheorem{lemma}[theorem]{Lemma}
\newtheorem{proposition}[theorem]{Proposition}

\theoremstyle{definition}

\theoremstyle{remark}
\newtheorem{remark}[theorem]{Remark}
\newtheorem*{acknowledgements}{Acknowledgements}

\numberwithin{equation}{section}

\newcommand{\acronym}[1]{\textsc{#1}}
\DeclareMathOperator{\Ad}{Ad}
\DeclareMathOperator{\ad}{ad}

\DeclareMathOperator{\CP}{\mathbb CP}

\DeclareMathOperator{\Gro}{\widetilde{Gr}}

\DeclareMathOperator{\Hol}{Hol}

\DeclareMathOperator{\Id}{Id}
\DeclareMathOperator{\im}{Im}
\newcommand{\inp}[3][]{\left\langle #2,#3\right\rangle_{#1}}

\DeclareMathOperator{\LieE}{\textsl E}
\DeclareMathOperator{\lieE}{\mathfrak e}
\DeclareMathOperator{\LieF}{\textsl F}

\DeclareMathOperator{\LieG}{\textsl G}
\DeclareMathOperator{\lieG}{\mathfrak g}
\DeclareMathOperator{\LSP}{\mathfrak{sp}}

\DeclareMathOperator{\Ort}{\textsl O}

\DeclareMathOperator{\re}{Re}

\DeclareMathOperator{\Sl}{\mathfrak{sl}}
\DeclareMathOperator{\SO}{\textsl{SO}}
\DeclareMathOperator{\so}{\mathfrak{so}}
\DeclareMathOperator{\SP}{\textsl{Sp}}

\newcommand{\subgroup}{\leqslant}
\DeclareMathOperator{\SU}{\textsl{SU}}
\DeclareMathOperator{\su}{\mathfrak{su}}

\DeclareMathOperator{\Un}{\textsl U}
\DeclareMathOperator{\un}{\mathfrak u}

\newenvironment{spmatrix}{\left(\smallmatrix}{\endsmallmatrix\right)}

\makeatletter

\newcommand{\ltwoa}{\mathord{\mathpalette\tw@a<}}
\newcommand{\rtwoa}{\mathord{\mathpalette\tw@a>}}
\newcommand{\tw@a}[2]{\ooalign{\hfil$#1 #2$\hfil\crcr
$#1 \Relbar\joinrel\Relbar$\crcr}}

\newcommand{\lthreea}{\mathord{\mathpalette\thr@a<}}
\newcommand{\rthreea}{\mathord{\mathpalette\thr@a>}}
\newcommand{\thr@a}[2]{\ooalign{\hfil$#1 #2$\hfil\crcr
\raise.3ex\hbox{$#1 \relbar\mathrel{\mkern-4mu}\relbar$}\crcr
$#1 \relbar\mathrel{\mkern-4mu}\relbar$\crcr
\lower.3ex\hbox{$#1 \relbar\mathrel{\mkern-4mu}\relbar$}\crcr}}

\newcommand{\Dyatop}[2]{\mathord{\mathpalette{\Dy@top{#1}{#2}}{}}}
\newcommand{\Dy@top}[3]{\overset{#3 #1}{#3 #2}}

\def\@captionfont{\normalfont\footnotesize}
\def\@captionheadfont{\scshape\footnotesize}

\makeatother

\begin{document}
\title[HyperKähler Potentials in Cohomogeneity Two]{HyperKähler Potentials
in\\ Cohomogeneity Two}

\author{Piotr Kobak}
\address[Kobak]{Insytut Matematyki\\
Uniwersytet Jagiello\'nski\\
ul.\ Reymonta 4\\
30-059 Kraków\\
Poland}
\email{kobak@im.uj.edu.pl}

\author{Andrew Swann}
\address[Swann]{Department of Mathematical Sciences\\
University of Bath\\
Claverton Down\\ 
Bath BA2 7AY\\
England}
\curraddr{Department of Mathematics and Computer Science\\
University of Southern Denmark\\
Odense University\\
Campusvej 55\\
DK-5230 Odense M\\
Denmark} 
\email{swann@imada.sdu.dk}


\begin{abstract}
  A hyperKähler potential is a function~$\rho$ that is a Kähler potential
  for each complex structure compatible with the hyperKähler structure.
  Nilpotent orbits in a complex simple Lie algebra are known to carry
  hyperKähler metrics admitting such potentials.  In this paper, we
  explicitly calculate the hyperKähler potential when the orbit is of
  cohomogeneity two.  In some cases, we find that this structure lies in a
  one-parameter family of hyperKähler metrics with Kähler potentials,
  generalising the Eguchi-Hanson metrics in dimension four.
\end{abstract}

\maketitle

\section{Introduction}
\label{sec:introduction}
HyperKähler metrics are special Ricci-flat structures that are known to
arise in many physical theories.  For example, moduli spaces of magnetic
monopoles often carry such metrics.  For good choices of boundary
conditions, these moduli spaces can be identified with more familiar
mathematical objects.  In this way, hyperKähler metrics have been shown to
exist on the adjoint orbits of a complex semi-simple Lie group~$G^{\mathbb
C}$ \cite{Kronheimer:nilpotent,Kronheimer:semi-simple,Biquard:orbits,
Kovalev:Nahm}.  In \cite{Dancer-Swann:hK-cohom1}, it was shown that these
examples include all hyperKähler metrics of cohomogeneity one.

Some of the earliest examples of hyperKähler metrics were found by
Calabi~\cite{Calabi:kaehler}.  His method was to take a complex symplectic
manifold, such as the cotangent bundle $T^*\CP(n)$, and find a potential
for a Kähler structure that would combine with the complex symplectic
structure to give a hyperKähler metric.  This approach has been applied to
certain semi-simple nilpotent orbits by a number of authors.  Biquard \&\ 
Gauduchon \cite{Biquard-Gauduchon:potential} gave a beautiful construction
for a potential on those semi-simple orbits that are the cotangent bundle
of a Hermitian symmetric space.  At the other extreme, Hitchin
\cite{Hitchin:integrable} used spectral theory to describe a potential for
the biggest semi-simple orbit in~$\Sl(n,\mathbb C)$ in terms of theta
functions (the special case of~$n=2$ may be found in~\cite{Santa-Cruz}).

Much attention has been paid to the semi-simple orbits, because one can
show that they are the only orbits to admit hyperKähler metrics that are
complete.  However, the incomplete metrics on nilpotent orbits still have
much interest.  One reason, is that each such orbit admits a
\emph{hyperKähler potential}, a function that is a Kähler potential for
each complex structure compatible with the hyperKähler structure, and so
these metrics on nilpotent orbits determine quaternionic Kähler metrics of
positive scalar curvature on a certain quotient manifold
\cite{Swann:MathAnn,Swann:HTwNil}.

The structures considered on coadjoint orbits are invariant under the action
of the compact group~$G$.  For nilpotent orbits, there is a natural partial
order given by inclusions of closures.  When $G$~is simple, the smallest
non-trivial orbits in this order are unique and they are distinguished by
being of cohomogeneity one under the action of~$G$.  In
\cite{Dancer-Swann:qK-cohom1}, it was shown that the nilpotent orbits of
cohomogeneity two also fit nicely in to the partial order: except when
$G=\SU(3)$, they are exactly the next-to-minimal orbits.  Given that the
nilpotent orbits of cohomogeneity one are understood
\cite{Dancer-Swann:hK-cohom1} (see also~\cite{Kobak-Swann:Wolf}), it is
natural to look at those of cohomogeneity two.  

In this paper, we consider cohomogeneity-two nilpotent orbits and find all
compatible $G$-invariant hyperKähler metrics on them that admit Kähler
potentials.  Our approach is that of Calabi's and we obtain the hyperKähler
potentials explicitly.  The hyperKähler potentials are unique, but in a few
cases we find that they lie in a one-parameter family of hyperKähler
metrics with Kähler potential.  These families may be regarded as
generalisations of the Eguchi-Hanson metrics in dimension four.

Combining our results with~\cite{Kobak-Swann:nilpotent}, means that
hyperKähler potentials are now known for all next-to-minimal orbits.  One
feature of the cohomogeneity-two case that makes the calculations possible,
is that each element of the orbit lies in a small rank~$2$ real subalgebra
which determines much of the hyperKähler structure.  In fact, unless $G$ is
the exceptional Lie group~$\LieG_2$, that subalgebra is $\so(4,\mathbb
C)=\Sl(2,\mathbb C)\oplus\Sl(2,\mathbb C)$ and the geometry is the product
of the structures from each factor.

For some cohomogeneity-two orbits the hyperKähler potential may also be
obtained by one of three other methods: a hyperKähler quotient
construction, a finite-cover by a minimal orbit for another group, or a
limit of a family of semi-simple orbits.  The first two methods will be
described elsewhere; the first only succeeds if the hyperKähler quotient is
sufficiently simple and the second only covers orbits on the list of
``shared orbits'' of Brylinski \&\ Kostant
\cite{Brylinski-Kostant:nilpotent-announcement}.  The third is contained in
Biquard \&\ Gauduchon's work~\cite{Biquard-Gauduchon:potential}.  However,
there are orbits for which the approach of this paper is the only one known
to give the result and our approach is uniform for all orbits of
cohomogeneity two.

\begin{acknowledgements}
  It is a pleasure to thank Brian Dupée for advice on \textsc{Maple},
  Alastair King and Francis Burstall for many useful conversations and
  Claude LeBrun for enlightenment.  This research was supported by the
  \acronym{epsrc} of Great Britain.  The first named author is also
  grateful for partial support from the KBN of Poland.
\end{acknowledgements}

\section{Preliminaries}
\label{sec:preliminaries}
\subsection{HyperKähler Structures}
\label{sec:hK}
Let $M$ be a manifold with endomorphisms $I$, $J$ and $K$ of the tangent
bundle~$TM$ satisfying the quaternion identities
\begin{equation*}
  I^2=J^2=-1\quad\text{and}\quad IJ=K=-JI.
\end{equation*}
This gives $T_xM$ the structure of an $\mathbb H$-module and so implies
that the dimension of~$M$ is a multiple of~$4$.  If $g$~is a Riemannian
metric on~$M$ preserved by $I$, $J$ and~$K$, in the sense that
$g(IX,IY)=g(X,Y)$, etc., for all tangent vectors $X,Y$, then we can define
two-forms $\omega_I$, $\omega_J$ and $\omega_K$ by
\begin{equation*}
  \omega_I(X,Y) = g(X,IY), \qquad\text{etc.}
\end{equation*}
If these three two-forms are closed, the structure $(M,g,I,J,K)$ is said to
be \emph{hyperKähler}.

Hitchin~\cite{Hitchin:Montreal} showed that on a hyperKähler manifold, the
almost complex structures $I$, $J$ and $K$ are integrable, and thus $(M,g)$
is a Kähler manifold in three distinct ways.  The restricted holonomy
group~$\Hol_g$ of~$(M,g)$ is then contained in~$\SP(n)$.  As $\SP(n)$ is a
subgroup of~$\SU(2n)$, this implies that any hyperKähler metric~$g$ is
Ricci-flat.

A function~$\rho\colon M\to \mathbb R$ is a Kähler potential for the
complex structure~$I$ if $\omega_I= -i\partial_I
\overline{\partial_I}\rho$.  This may be reformulated as
\begin{equation}
  \label{eq:dId}
  \begin{split}
    \omega_I
    &= - i \partial_I \overline{\partial_I} \rho
    = - i d \overline{\partial_I} \rho = - \tfrac i2 d (d-iId) \rho \\
    &= - \tfrac 12 d I d \rho.
\end{split}
\end{equation}
The function~$\rho$ is a \emph{hyperKähler potential} if it is
simultaneously a Kähler potential for $I$, $J$ and~$K$.  HyperKähler
potentials are defined up to an additive constant.  The existence of a
hyperKähler potential implies strong restrictions on the geometry of~$M$
\cite{Swann:MathAnn}: the metric~$g$ and potential~$\rho$ satisfy
$\nabla^2\rho=g$; the manifold~$M$ admits an infinitesimal action of
$\mathbb H^*$, with $\SP(1)\subgroup\mathbb H^*$ preserving~$g$ and
permuting $I$, $J$ and $K$; the $\mathbb H^*$-orbits are flat and totally
geodesic; locally $M$ fibres over a quaternionic Kähler orbifold of
positive scalar curvature.

We will be considering hyperKähler structures that are invariant under the
action of a compact group~$G$.  It is therefore worth noting that if we
have a Kähler potential then this may be taken to be $G$-invariant.
Indeed, if $\rho$~is any Kähler potential, then since the $G$-action
preserves~$I$, the expression $\partial_I \overline{\partial_I} \rho$ is
equivariant for the action of~$G$.  However, $\omega_I= - i \partial_I
\overline{\partial_I} \rho$, is assumed to be $G$-invariant, so averaging
$\rho$ over the group action produces an invariant Kähler potential.

\subsection{Lie Algebras and Orbits}
\label{sec:Lie}
On the semi-simple complex Lie algebra~$\mathfrak g^{\mathbb C}$, let
$\inp\cdot\cdot=\inp[\mathfrak g]\cdot\cdot$ be the negative of the Killing
form and let $\sigma$~be a real structure giving a compact real
form~$\mathfrak g$ of~$\mathfrak g^{\mathbb C}$.

At a point~$X$ of a nilpotent orbit~$\mathcal O$, the vector field
generated by~$A$ in~$\mathfrak g^{\mathbb C}$ is $\xi_A=[A,X]$.  These
vector fields satisfy $[\xi_A,\xi_B] = \xi_{-[A,B]}$, for all
$A,B\in\mathfrak g^{\mathbb C}$.

The orbit~$\mathcal O$ carries a complex structure~$I$ defined by
\begin{equation}
  \label{eq:I}
  I\xi_A = i \xi_A = \xi_{iA}.
\end{equation}
There is also a complex symplectic form, known as the
Kirillov-Kostant-Souriau form, on~$\mathcal O$ which we take to be given by
\begin{equation}
  \label{eq:omega-c}
  \omega_c^{\mathcal O}(\xi_A,\xi_B)_X = \inp X{[A,B]} = - \inp{\xi_A}B.
\end{equation}
We will be looking for hyperKähler structures on~$\mathcal O$ with $I$
given by~\eqref{eq:I} and $\omega_J +i\omega_K=\omega_c^{\mathcal O}$.  We
will call these \emph{compatible hyperKähler structures} on~$\mathcal O$.

\section{Potentials Depending on Two Invariants}
\label{sec:2-i}
Consider the following two functions on a nilpotent orbit~$\mathcal O$:
\begin{equation*}
  \eta_1(X) = \inp X{\sigma X} \quad\text{and}\quad \eta_2(X) = - \inp
  {[X,\sigma X]}{[X,\sigma X]}.
\end{equation*}
Note that $\eta_2(X)=\inp Y{\sigma Y}$ with $Y=[X,\sigma X]$, so is
positive, and that both $\eta_1$ and $\eta_2$ are invariant under the
action of the compact group~$G$.  Suppose $\rho$ is a Kähler potential
for~$I$ depending only on $\eta_1$ and~$\eta_2$, i.e.,
\begin{equation}
  \label{eq:rho}
  \rho = \rho (\eta_1,\eta_2).
\end{equation}

\begin{lemma}
  At $X\in\mathcal O$, the two-form $\omega_I$ defined by~$\rho$
  in formula \eqref{eq:dId} is
  \begin{equation}
    \label{eq:omega-I}
    \begin{split}
      \omega_I(\xi_A,\xi_B)_X
      &= 2 \rho_1 \im \inp {\xi_A}{\sigma\xi_B} \\
      &\quad - 4 \rho_2 \im \inp {\xi_A}{[\sigma\xi_B,[X,\sigma X]]
      + [\sigma X,[X,\sigma\xi_B]]} \\
      &\quad + 2 \rho_{11} \im
      \bigl(
      \inp {\xi_A}{\sigma X}\inp {\sigma\xi_B}X
      \bigr)
      \\
      &\quad  - 4 \rho_{12}
      \begin{aligned}[t]
        \im
        \bigl( &
        \inp {\xi_A}{[\sigma X,[X,\sigma X]]}\inp{\sigma \xi_B}X
        \\
        &+ \inp {\xi_A}{\sigma X} \inp {\sigma \xi_B}{[X,[\sigma X,X]]}
        \bigr)
      \end{aligned}
      \\
      &\quad + 8 \rho_{22} \im
      \bigl(
      \inp {\xi_A}{[\sigma X,[X,\sigma X]]}
      \inp {\sigma \xi_B}{[X,[\sigma X,X]]}
      \bigr)
      ,
    \end{split}
  \end{equation}
  where $\rho_i = \partial \rho / \partial \eta_i$, etc.
\end{lemma}

\begin{proof}
  Expanding~\eqref{eq:dId}, we have
  \begin{multline}
    \label{eq:2d}
    -2\omega_I 
    = \rho_1\, dId\eta_1 + \rho_2\, dId\eta_2
    + \rho_{11}\, d\eta_1\wedge Id\eta_1 \\
    + \rho_{12} (d\eta_2\wedge Id\eta_1 + d\eta_1 \wedge Id\eta_2)
    + \rho_{22} d\eta_2 \wedge Id\eta_2.
  \end{multline}
  The exterior derivative of~$\eta_1$ is given by
  \begin{equation*}
    d\eta_1(\xi_A)_X
    = \inp{[A,X]}{\sigma X} + \inp X{\sigma[A,X]}
    = 2 \re \inp{\xi_A}{\sigma X}.
  \end{equation*}
  Hence $Id\eta_1(\xi_A) = 2\im\inp{\xi_A}{\sigma X}$ and
  $dId\eta(\xi_A,\xi_B) = - 4 \im \inp{\xi_A}{\sigma \xi_B}$,
  at~$X\in\mathcal O$.
  
  For $\eta_2$, the initial computation is similar and gives
  \begin{equation*}
    d\eta_2(\xi_A)_X = - 4\re \inp{\xi_A}{[\sigma X,[X,\sigma X]]}.
  \end{equation*}
  The second derivative, however, is slightly more involved:
  \begin{equation*}
    \begin{split}
      dId&\eta_2(\xi_A,\xi_B)_X \\
      & = \xi_A(Id\eta_2(\xi_B)) - \xi_B(Id\eta_2(\xi_A)) -
      Id\eta_2([\xi_A,\xi_B]) \\
      & = -4 
      \begin{aligned}[t]
        \im
        \bigl\{
        &\inp{\xi_B}{[\sigma \xi_A,[X,\sigma X]]}
        + \inp{\xi_B}{[\sigma X,[\xi_A,\sigma X]]}\\
        &\quad+ \inp{\xi_B}{[\sigma X,[X,\sigma \xi_A]]}
        + \inp{[B,\xi_A]}{[\sigma X,[X,\sigma X]]}\\
        & -\inp{\xi_A}{[\sigma \xi_B,[X,\sigma X]]}
        - \inp{\xi_A}{[\sigma X,[\xi_B,\sigma X]]}\\
        &\quad- \inp{\xi_A}{[\sigma X,[X,\sigma \xi_B]]}
        - \inp{[A,\xi_B]}{[\sigma X,[X,\sigma X]]}\\
        & + \inp{[[A,B],X]}{[\sigma X,[X,\sigma X]]}
        \bigr\}
      \end{aligned}
      \\
      & = -4
      \begin{aligned}[t]
        \im
        \bigl\{
        & - \inp{[\xi_A,\sigma \xi_B]}{[X,\sigma X]}
        + \inp{[\xi_B,\sigma \xi_A]}{[X,\sigma X]} \\
        & + \inp{[\sigma X, \xi_A]}{[X,\sigma \xi_B]}
        - \inp{[X,\sigma \xi_A]}{[\sigma X,\xi_B]}
        \bigr\}
      \end{aligned}
      \\
      & = 8 \im \inp{\xi_A}{[\sigma \xi_B,[X,\sigma X]] + [\sigma
      X,[X,\sigma \xi_B]]}.
    \end{split}
  \end{equation*}
  Combining these formulæ gives the claimed result.
\end{proof}

The two-form $\omega_I$ is our candidate for a Kähler form on~$\mathcal O$.

\begin{remark}
  \label{rem:metric}
  The corresponding symmetric bilinear form is given by
  $g(\xi_A,\xi_B)=\omega_I(I\xi_A,\xi_B)$ and is simply the right-hand side
  of equation~\eqref{eq:omega-I} with `$\im$' replaced by `$\re$'
  throughout.
\end{remark}

We will eventually require~$g$ to be positive definite.  However for now
simply assume that $g$~is non-degenerate and define an endomorphism~$J$
of~$T_X\mathcal O$ by
\begin{equation}
  \label{eq:g-J}
  g(\xi_A,\xi_B) = \re \omega_c^{\mathcal O}(J\xi_A,\xi_B).
\end{equation}

\begin{lemma}
  The endomorphism~$J$ of~$T_X\mathcal O$ is given by
  \begin{equation}
    \label{eq:J}
    \begin{split}
      J\xi_A
      & = -2 \rho_1 [X,\sigma \xi_A] \\
      & \quad+ 4 \rho_2 (2 [X,[\sigma X,[X,\sigma \xi_A]]] - [X,[X,[\sigma
      X,\sigma \xi_A]]])\\
      & \quad - 2\rho_{11} \inp{\sigma \xi_A}X [X,\sigma X]\\
      & \quad + 4\rho_{12}
      \begin{aligned}[t]
        \bigl(
        & \inp{\sigma \xi_A}{[X,[\sigma X,X]]} [X,\sigma X] \\
        &+ \inp{\sigma \xi_A}X [X,[\sigma X,[X,\sigma X]]]
        \bigr)
      \end{aligned}
      \\
      & \quad - 8 \rho_{22} \inp{\sigma \xi_A}{[X,[\sigma X,X]]} [X,[\sigma
      X,[X,\sigma X]]].
    \end{split}
  \end{equation}
\end{lemma}

\begin{proof}
  Equation~\eqref{eq:g-J} implies $g(\xi_A,\xi_B)=-\re\inp{J\xi_A}B$, and
  then \eqref{eq:omega-I} gives the above formula for~$J$, except that the
  coefficient of~$\rho_2$ is
  \begin{equation}
    \label{eq:rho2-coefficient}
    4\bigl( [X,[\sigma \xi_A,[X,\sigma X]]] + [X,[\sigma X,[X,\sigma
    \xi_A]]] \bigr).
  \end{equation}
  Using the Jacobi identity, we have
  \begin{equation*}
    [\sigma \xi_A,[X,\sigma X]]
    = - [X,[\sigma X,\sigma\xi_A]] + [\sigma X,[X,\sigma\xi_A]].
  \end{equation*}
  Applying this to the first term in~\eqref{eq:rho2-coefficient} gives the
  result.
\end{proof}

At this stage there is no guarantee that $J^2=-1$.  It is imposing this
condition that severely restricts the possibilities for~$\rho$.

\section{Small Nilpotent Orbits and Real Subalgebras}
\label{sec:subalgebras}
The nilpotent orbits in~$\mathfrak g^{\mathbb C}$ are partially ordered by
saying $\mathcal O_1 \preceq \mathcal O_2$ if and only if $\mathcal O_1
\subset \overline{\mathcal O_2}$.  When $\mathfrak g^{\mathbb C}$~is
simple, there is a unique non-zero orbit~$\mathcal O_{\text{min}}$ which is
minimal for this partial order.  This orbit is of cohomogeneity one with
respect to the action of the compact group~$G$, and for each $X\in\mathcal
O_{\text{min}}$, the subalgebra spanned by $\{X,\sigma X\}$ is isomorphic
to $\Sl(2,\mathbb C)$ and is the complexification of an~$\su(2)$-subalgebra
of~$\mathfrak g$.  Note that $\mathcal O_{\text{min}}$ is the orbit of a
root vector for the longest root.

In general, the Jacobsen-Morosov Theorem says that each nilpotent
element~$X$ lies in an $\Sl(2,\mathbb C)$-subalgebra (see
e.g.~\cite{Carter:finite}).  However, in general this subalgebra is not
$\sigma$-invariant.  The following result is usually attributed to
Borel~\cite{Borel:proper}.

\begin{proposition}[Borel]
  \label{prop:real-su2}
  Each nilpotent orbit~$\mathcal O$ contains an element $X$ such that the
  linear span of $\{X,\sigma X, [X,\sigma X]\}$ is a real subalgebra
  isomorphic to~$\Sl(2,\mathbb C)$.
\end{proposition}

\begin{proof}
  Fix $X'$ in $\mathcal O$ and take any $\Sl(2,\mathbb C)$ containing~$X$.
  There are $H$ and~$Y$ in~$\Sl(2,\mathbb C)$ such that $[H,X']=2X'$,
  $[X',Y]=H$ and $[H,Y]=-2Y$.  The element~$H$ is thus semi-simple
  in~$\Sl(2,\mathbb C)$ and hence in $\mathfrak g^{\mathbb C}$, so we find
  a Cartan subalgebra~$\mathfrak t$ of~$\mathfrak g^{\mathbb C}$
  containing~$H$ and choose a system of positive roots~$\Delta^+$ so that
  $X$ lies in a sum of positive root spaces.  The pair $(\mathfrak
  t,\Delta^+)$ has an associated real structure~$\sigma'$, which maps
  $\Delta^+$ to~$\Delta^-$ and defines a compact real form of~$\mathfrak
  g^{\mathbb C}$.  Now all compact real forms of~$\mathfrak g^{\mathbb C}$
  are conjugate, so there is a $g\in G^{\mathbb C}$ such that
  $\Ad_g(\sigma'A)=\sigma\Ad_g A$, for all $A\in\mathfrak g^{\mathbb C}$.
  Taking $X=\Ad_g X'$ gives an element of~$\mathcal O$ of the desired type.
\end{proof}

Let us recall the Morse theory picture of the nilpotent variety described
in~\cite{Swann:HTwNil} (see also \cite{Kobak-Swann:nilpotent,
Kronheimer:nilpotent}).  Each nilpotent orbit~$\mathcal O$ admits a certain
free action of $\mathbb H^*/\{\pm1\}$.  The quotient $\mathfrak M(\mathcal
O)=\mathcal O/\mathbb H^*$ may be described as a submanifold of the
Grassmannian $\Gro_3(\mathfrak g)$ of oriented three-planes in the real Lie
algebra~$\mathfrak g$.  One defines a functional
$\psi\colon\Gro_3(\mathfrak g)\to\mathbb R$ by
$\psi(V)=\inp{e_1}{[e_2,e_3]}$, where $\{e_1,e_2,e_3\}$ is an oriented
orthonormal basis for~$V$.  Away from zero, $\psi$~is a non-degenerate
$G$-equivariant Morse function in the sense of Bott.  The points on the
non-zero critical sets correspond to subalgebras of~$\mathfrak g$
isomorphic to~$\su(2)$.  The set of real $\su(2)$-subalgebras associated
to~$\mathcal O$ via Proposition~\ref{prop:real-su2}, oriented so that
$\psi$ is positive, forms a non-zero critical manifold~$\mathcal C(\mathcal
O)$.  The manifold $\mathfrak M(\mathcal O)$ is the stable manifold
attached to~$\mathcal C(\mathcal O)$.  The partial order on stable
manifolds for the gradient flow induces the partial order~$\preceq$ on
nilpotent orbits.  In particular, the maximum of~$\psi$ is achieved
on~$\mathfrak M(\mathcal O_{\text{min}})$.

We are interested in orbits of cohomogeneity two.  These were computed
in~\cite{Dancer-Swann:qK-cohom1} and are the orbits listed in
Table~\ref{tab:c-2}.  We say that a nilpotent orbit~$\mathcal O$ is
\emph{next-to-minimal} if $\mathcal O \succneqq \mathcal O_{\text{min}}$
and there is no orbit~$\mathcal O'$ with $\mathcal O \succneqq \mathcal O'
\succneqq \mathcal O_{\text{min}}$.  It is pleasing to note that the orbits
listed in Table~\ref{tab:c-2} are precisely the next-to-minimal orbits in
the given algebras.  The only next-to-minimal orbit that does not occur is
that in~$\Sl(3,\mathbb C)$, which is cohomogeneity four.

\begin{table}[tbp]
  \begin{center}
    \newcommand{\tablestrut}{\vrule height 12pt width 0pt}
    \leavevmode
    \begin{tabular}[t]{lc}
      \hline
      Type & Orbit \\
      \hline \hline
      \tablestrut $A_n$ & $(2^21^{n-3})$ 
      \\
      \hline
      \tablestrut \rlap{$B_{(n-1)/2}$, $D_{n/2}$} &\\
      & $(31^{n-3})$ 
      \\
      & $(2^41^{n-8})$  \\
      \hline
      \tablestrut $C_n$ & $(2^21^{2n-4})$  \\ 
      \hline
    \end{tabular}
    \qquad
    \begin{tabular}[t]{lc}
      \hline
      Type & Orbit \\
      \hline \hline
      \tablestrut $\LieG_2$ & $\scriptstyle 0\rthreea1$ \\
      \tablestrut $\LieF_4$ & $\scriptstyle 00\rtwoa01$  \\
      \tablestrut $\LieE_6$ & $\scriptstyle 10\Dyatop0001$ \\ 
      \tablestrut $\LieE_7$ & $\scriptstyle 010\Dyatop0000$ \\
      \tablestrut $\LieE_8$ & $\scriptstyle 0000\Dyatop0001$ \\
      \hline
    \end{tabular}
    \caption{Orbits of cohomogeneity two in simple Lie algebras.  Orbits in
    classical algebras are specified by partitions and $n$~is to be taken
    large enough so that the partition can occur.  The orbits in
    exceptional algebras are given by their weighted Dynkin diagram (see
    e.g.~\cite{Carter:finite}).  Note that for type $D_{2m}$, the partition
    $(2^41^{4m-8})$ describes two orbits; their union is one orbit under
    the action of~$\Ort(2m)$.}
    \label{tab:c-2}
  \end{center}
\end{table}

Recall that according to Proposition~\ref{prop:sl2-potentials} elements of
cohomogeneity-one nilpotent orbits lie in a real $\Sl(2,\mathbb C)$, i.e.,
in a $\sigma$-invariant rank one Lie algebra.  It is remarkable that the
elements of cohomogeneity-two orbits lie in $\sigma$-invariant rank two Lie
algebras.  The following can be thought of as a cohomogeneity-two version
of Borel's result.

\begin{theorem}
  \label{thm:so4}
  Suppose $G$ is a compact simple Lie group and that $\mathcal O$ is a
  nilpotent orbit in~$\mathfrak g^{\mathbb C}$ of cohomogeneity two.
  Suppose $X$ is an element of~$\mathcal O$ that does not lie in a real
  $\Sl(2,\mathbb C)$-subalgebra.
  
  Let $\mathfrak h_X^{\mathbb C}$ be the subalgebra of~$\mathfrak
  g^{\mathbb C}$ generated by $X$ and~$\sigma X$.  Then $\mathfrak
  h_X^{\mathbb C}$ is isomorphic to~$\so(4,\mathbb C)$, unless $\mathfrak
  g=\lieG_2$, in which case $\mathfrak h^{\mathbb C}\cong\lieG_2^{\mathbb
  C}$.  In all cases, the embedding $\mathfrak h^{\mathbb C}
  \hookrightarrow \mathfrak g^{\mathbb C}$ is a homothety with respect to
  the Killing forms.
\end{theorem}

\begin{proof}
  Consider the Morse theory picture.  Firstly, in $\mathfrak g^{\mathbb
  C}$, the closure of~$\mathcal O$ is $\mathcal O\cup \mathcal
  O_{\text{min}} \cup \{0\}$.  In $\Gro_3(\mathfrak g)$ we have
  $\overline{\mathfrak M(\mathcal O)}=\mathfrak M(\mathcal O)\cup\mathfrak
  M(\mathcal O_{\text{min}})$.  For the orbits of cohomogeneity two,
  $\mathfrak M(\mathcal O)$~is a manifold of cohomogeneity one; the usual
  scaling by $\mathbb R_{>0}$, which is also part of the $\mathbb
  H^*/\{\pm1\}$-action, is transverse to the $G$-orbits on~$\mathcal O$.
  Suppose $\ell$~is a curve in $\overline{\mathfrak M(\mathcal O)}$ joining
  a point of~$\mathcal C(\mathcal O)$ to a point of~$\mathfrak M(\mathcal
  O_{\text{min}})$.  Then the fact that $\mathfrak M(\mathcal O)$ is the
  stable manifold for the gradient flow of the $G$-invariant
  functional~$\psi$, implies that the image of~$\ell$ in $\mathfrak
  M(\mathcal O)/G$ is the whole (one-dimensional) quotient space.
  
  Now to parameterise $\overline{\mathcal O}/G$ it is enough to find a
  two-dimensional family of elements which is invariant under scaling and
  contains an element lying over $\mathcal C(\mathcal O)$ and an element
  of~$\mathcal O_{\text{min}}$.  When $\mathfrak g\ne\lieG_2$, we will find
  such a family lying in a $\sigma$-invariant $\so(4,\mathbb C)$-subalgebra
  of~$\mathfrak g^{\mathbb C}$.
  
  The Lie algebra $\so(4,\mathbb C)$ splits as $\Sl(2,\mathbb
  C)_+\oplus\Sl(2,\mathbb C)_-$.  It contains three non-trivial nilpotent
  orbits: $\mathcal O_\pm$, the non-trivial nilpotent orbits in the
  factor~$\Sl(2,\mathbb C)_\pm$; and $\mathcal O_\Delta=\mathcal
  O_+\times\mathcal O_-$.  The orbits $\mathcal O_{\pm}$ are cohomogeneity
  one and $\mathcal O_\Delta$ is cohomogeneity two.  Our orbit~$\mathcal O$
  will meet~$\so(4,\mathbb C)$ in~$\mathcal O_\Delta$ and $\mathcal
  O_+\cup\mathcal O_-$ will be the intersection $\mathcal
  O_{\text{min}}\cap\so(4,\mathbb C)$.

  For the classical groups, we use the Jordan normal forms for elements of
  the orbits.  For type~$A_n$,  the Jordan normal form
  is~$(2^21^{n-3})$ and the matrices
  \begin{equation*}
    X_{s,t}=
    \begin{spmatrix}
      0&s&&& \\
      0&0&&& \\
      &&0&t& \\
      &&0&0& \\
      &&&&\text{\large 0}
    \end{spmatrix}
  \end{equation*}
  lie in the orbit unless $s$ or $t$ is zero.  They also lie in the
  $\su(2)\oplus\su(2)$-subalgebra contained in the first two $(2\times2)$
  diagonal blocks.  The matrix $X_{1,1}$ lies in a real $\Sl(2,\mathbb
  C)$-subalgebra and $X_{1,0}$ is in $\mathcal O_{\text{min}}$.  So this
  two parameter family is as required.  Exactly the same technique works
  for $C_n$.
  
  For types $B$ and~$D$, we are looking at matrices in~$\so(n,\mathbb C)$.
  It is convenient to take $\so(n,\mathbb C)$ to be the set of complex
  $(n\times n)$ matrices~$A$ such that $A^tB+BA=0$, where $B$~is the matrix
  with $1$'s down the anti-diagonal and $0$'s elsewhere.  For Jordan form
  $(31^{n-3})$, we just take an $\so(4,\mathbb C)$-subalgebra containing
  the Jordan block~$(3)$.  When the Jordan type is $(2^41^{n-8})$, and the
  Lie algebra type is not~$D_{2n}$, we have the same situation as
  for~$A_n$, but now the blocks come in pairs.  Thus the two families one
  considers are
  \begin{equation}
    \label{eq:B-D}
    \begin{spmatrix}
      0&s&&&&&&& \\
      0&0&&&&&&& \\
      &&0&t&&&&& \\
      &&0&0&&&&& \\
      &&&&\text{\large 0}&&&&\\
      &&&&&0&-t&& \\
      &&&&&0&0&& \\
      &&&&&&&0&-s \\
      &&&&&&&0&0 
    \end{spmatrix}
    \quad
    \text{and}
    \quad
    \begin{spmatrix}
      \text{\large 0}&&&&&\text{\large 0}\\
      &0&s&t&0& \\
      &&0&0&-t& \\
      &&&0&-s& \\
      &&&&0& \\
      \text{\large 0}&&&&&\text{\large 0}
    \end{spmatrix}
    .
  \end{equation}
  
  For $D_{2n}$, the matrices of Jordan type $(2^41^{n-8})$ form a single
  $\Ort(n,\mathbb C)$-orbit but split into two orbits~$(2^41^{n-8})_\pm$
  under the action of~$\SO(n,\mathbb C)$.  We thus obtain $(2^41^{n-8})_-$
  from $(2^41^{n-8})_+$ by conjugating by an element~$W$ of
  determinant~$-1$ in~$\Ort(n,\mathbb C)$.  One now considers three
  representative matrices, two as in~\eqref{eq:B-D} and
  \begin{equation*}
      \begin{spmatrix}
      0&0&&&&&&& \\
      0&0&&&&&&& \\
      &&0&t&&&&& \\
      &&0&0&&&&& \\
      &&&&\text{\large 0}&&&&\\
      &&&&&0&-t&& \\
      &&&&&0&0&& \\
      s&0&&&&&&0&0 \\
      0&-s&&&&&&0&0 
      \end{spmatrix}
    , \quad\text{obtained with }
    W =
      \begin{spmatrix}
        0&&&&-1\\
        &1&&&\\
        &&\ddots&&\\
        &&&1&\\
        -1&&&&0
      \end{spmatrix}
      .
  \end{equation*}
  
  In all cases, the matrix lies over $\mathcal C(\mathcal O)$ when $s=t\ne
  0$ and is in $\mathcal O_{\text{min}}$ when $t=0$ and $s\ne 0$.
  
  For the exceptional Lie algebras we use the Beauville bundle~$N(\mathcal
  O)$ \cite{Beauville:Fano} as a tool for computation.  This bundle is
  defined as follows.  Find a real $\Sl(2,\mathbb C)$-subalgebra associated
  to~$\mathcal O$ and let $\{e,f,h\}$ be a basis for this subalgebra, with
  $f=-\sigma e$, $h=[e,f]$ and $[h,e]=2e$.  The eigenvalues of~$\ad h$
  on~$\mathfrak g^{\mathbb C}$ are known to be integers
  (see~\cite{Carter:finite}).  Let $\mathfrak g(i)$ be the $i$-eigenspace
  of~$\ad h$.  Put
  \begin{equation*}
    \mathfrak p = \bigoplus_{i\geqslant 0} \mathfrak g(i)
    \quad\text{and}\quad
    \mathfrak n = \bigoplus_{i\geqslant 2} \mathfrak g(i).
  \end{equation*}
  Then $\mathfrak p$ is a parabolic subalgebra of~$\mathfrak g^{\mathbb C}$
  and the corresponding homogeneous space $\mathcal F=G^{\mathbb C}/P$ is a
  flag manifold.  The subalgebra~$\mathfrak n$ is preserved by the adjoint
  action of~$P$ and the Beauville bundle~$N(\mathcal O)$ is defined to be
  the bundle over~$\mathcal F$ associated to~$\mathfrak n$, i.e.,
  \begin{equation*}
    N(\mathcal O) = G^{\mathbb C}\times_P \mathfrak n.
  \end{equation*}
  The important property of~$N(\mathcal O)$ is that it contains~$\mathcal
  O$ as an open dense $G^{\mathbb C}$-orbit.
  
  Now each flag manifold is a homogeneous manifold for the action of the
  compact group.  So $\mathcal F= G/K$ for some compact subgroup~$K$
  of~$G$.  (In fact, the Lie algebra~$\mathfrak k$ of~$K$ is given by
  $\mathfrak k^{\mathbb C}=\mathfrak g(0)$.)  The Beauville bundle is then
  $G\times_K\mathfrak n$ and the cohomogeneity of~$\mathcal O$ is the
  cohomogeneity of the action of~$K$ on~$\mathfrak n$.  Choose a Cartan
  subalgebra in~$\mathfrak g(0)$ and a root system for~$\mathfrak g(0)$
  with all root spaces in~$\mathfrak p$.  Note that, by definition, the
  weighted Dynkin diagram for~$\mathcal O$ gives the eigenvalues of~$\ad h$
  on the positive simple root spaces, from which all the other eigenvalues
  are easily computed.
  
  In the case of cohomogeneity-two orbits not in~$\lieG_2$, we find that
  $\mathfrak n\cong \mathbb R^2\otimes V$ as a representation
  of~$K=\SO(2)L$, with $V$~irreducible and $L$~acting two-point
  transitively on the unit sphere in~$V$.  Thus under the action of~$K$, we
  can move any nilpotent element in~$\mathcal O$ into any complex
  two-dimensional subspace (for the complex structure induced by the action
  of~$\SO(2)$).  We then find root spaces $\mathfrak g_\alpha$ and
  $\mathfrak g_\beta$ contained in~$\mathfrak n$ with $\alpha$ and $\beta$
  orthogonal long roots such that $\alpha\pm\beta$ is not a root.  The
  $\sigma$-invariant subalgebra containing these root spaces is then the
  required~$\so(4,\mathbb C)$.  For the relevant four exceptional algebras,
  this information is given in Table~\ref{tab:so4}.
  
  \begin{table}[tbp]
    \begin{center}
      \newcommand{\tablestrut}{\vrule height 12pt width 0pt}
      \begin{tabular}{ccccc}
        \hline
        Type & $\mathfrak k$&$\mathfrak n$&$\alpha$&$\beta$\\
        \hline \hline
        \tablestrut $F_4$&$\so(2)+\so(7)$&$\mathbb R^2\otimes\mathbb R^7$
        &$\scriptstyle 23\rtwoa42$&$\scriptstyle 01\rtwoa22$\\
        \tablestrut $E_6$&$2\so(2)+\so(8)$&$\mathbb R^2\otimes\mathbb R^8$
        &$\scriptstyle 12\Dyatop2321$&$\scriptstyle 11\Dyatop1221$\\
        \tablestrut $E_7$&$\so(2)+\su(2)+\so(10)$&$\mathbb
        R^2\otimes\mathbb R^{10}$ 
        &$\scriptstyle 123\Dyatop2432$&$\scriptstyle 122\Dyatop2321$\\
        \tablestrut $E_8$&$\so(2)+\so(14)$&$\mathbb R^2\otimes\mathbb
        R^{14}$&$\scriptstyle 2345\Dyatop3642$&$\scriptstyle 0123\Dyatop
        2432$\\ 
        \hline
      \end{tabular}
      \caption{The data for $\so(4)$-subalgebras corresponding to
      next-to-minimal orbits in four exceptional algebras.} 
      \label{tab:so4}
    \end{center}
  \end{table}

  For $G_2$, the isotropy group for~$\mathcal F$ is $K=\Un(1)\SU(2)$.  We
  have $\mathfrak n=\mathfrak g(2)+\mathfrak g(3)$ with $\mathfrak g(2)\cong
  L^2$ and $\mathfrak g(3)\cong L^3 S^1$, where $L=\mathbb C$ and
  $S^1=\mathbb C^2$ are the fundamental representations of~$\Un(1)$ and
  $\SU(2)$, respectively.  The orbit $\mathcal O$ in this case is the orbit
  of short root vectors and the subalgebra generated by $X$ and $\sigma X$
  contains both short and long roots, so is all of~$\lieG_2$.
\end{proof}

\section{Two Models}
\label{sec:models}
In this section we compute Kähler potentials for hyperKähler structures on
two particular nilpotent orbits: one in~$\Sl(2,\mathbb C)$ and the other
in~$\so(4,\mathbb C)$.  These results will be used in the next section to
derive the hyperKähler potentials for cohomogeneity-two orbits.  In view of
Theorem~\ref{thm:so4}, we consider these cases with inner products that are
multiples of that given by the Killing form.

We start by considering $\mathfrak g^{\mathbb C}=\Sl(2,\mathbb C)$ with
inner product $k^2\inp[\Sl(2)]\cdot\cdot$, where $k>0$~is constant and
$\inp[\Sl(2)]\cdot\cdot$ is negative of the Killing form.  This Lie algebra
contains only one non-trivial nilpotent orbit~$\mathcal O$ consisting of
the $(2\times2)$ matrices~$X$ such that $X^2=0$ and $X\ne0$.  The orbit is
the minimal nilpotent orbit in~$\Sl(2,\mathbb C)$ and is of cohomogeneity
one under the adjoint action of~$\SU(2)$.  In fact two elements of the
orbit have the same norm if and only if they are $\SU(2)$-conjugate.  Thus
any $\SU(2)$-invariant Kähler potential~$\rho$ on~$\mathcal O$ is a
function of just~$\eta=k^2\inp[\Sl(2)] X{\sigma X}$.  Write
\begin{equation}
  \label{eq:efh}
  e = 
  \begin{pmatrix}
    0&1\\
    0&0
  \end{pmatrix}
  ,\quad
  f =
  \begin{pmatrix}
    0&0\\
    1&0
  \end{pmatrix}
  \quad\text{and}\quad
  h = 
  \begin{pmatrix}
    1&0\\
    0&-1
  \end{pmatrix}
  .
\end{equation}
Then $\{e,f,h\}$ is an~$\Sl(2,\mathbb C)$ triple, with $f=-\sigma e$
and~$h=-\sigma h$.  Using the action of~$\SU(2)$, we may assume that
$X=t\,e$, for some $t>0$.  The tangent space $T_X\mathcal
O=[X,\Sl(2,\mathbb C)]$ is spanned by $e$ and~$h$.  If we consider the
complex symplectic form $k^2\omega_C^{\mathcal O}$, we may perform the same
calculations as in~\eqref{eq:J} and get
\begin{equation*}
  J_X e = 2 t (\rho' + \eta \rho'') h
  \quad\text{and}\quad
  J_X h = -4 t \rho' e ,
\end{equation*}
where $\rho'=d\rho/d\eta$, etc.  As $\eta(X) = 4k^2t^2$, the condition that
$J^2=-1$ is equivalent to
\begin{equation}
  \label{eq:ode}
  2 \eta \rho' (\rho' + \eta \rho'') = k^2.
\end{equation}
The left-hand side is simply the derivative of $(\eta\rho')^2$ with respect
to~$\eta$, so
\begin{equation}
  \label{eq:rhop2-sl2}
  {\rho'}^2=(k^2\eta+c)/\eta^2,
\end{equation}
for some constant~$c$.  In order to have the potential defined on the whole
orbit we need $c\geqslant 0$.  The corresponding metric may be calculated
as in Remark~\ref{rem:metric} and is given by
\begin{equation}
  \label{eq:metric}
  \begin{split}
    g(\xi_A,\xi_B) = \frac{2k^4}{\eta} \re
    \Bigl(&\rho'\bigl(\inp{\xi_A}{\sigma\xi_B}\inp X{\sigma X}-
    \inp{\xi_A}{\sigma X}\inp{X}{\sigma \xi_B}\bigr)\\
    &\quad +\frac{k^2}{2\eta\rho'}\inp{\xi_A}{\sigma X}\inp{X}{\sigma
    \xi_B}\Bigr),
\end{split}
\end{equation}
which is positive definite provided we take the positive square root
in~\eqref{eq:rhop2-sl2}.  Now \eqref{eq:rhop2-sl2} determines $\rho$ up to
an additive constant, and this is enough to fix the metric structure.

\begin{proposition}
  \label{prop:sl2-potentials}
  For fixed~$k$, the nilpotent orbit in~$\Sl(2,\mathbb C)$ has a
  one-parameter family of $\SU(2)$-invariant hyperKähler metrics with a
  Kähler potential and with $k^2\omega_c^{\mathcal O}$ as the complex
  symplectic form.  The $G$-invariant Kähler potential $\rho$ is given by
  \begin{equation}
    \label{eq:rhop-sl2}
    \rho' = \frac 1\eta \sqrt{k^2\eta+c},
  \end{equation}
  where $c\geqslant0$ is a constant and $\eta(X)=k^2\inp[\Sl(2)] X{\sigma
  X}$.  \qed
\end{proposition}

\noindent
Note that if we rewrite everything in terms of the variable~$t$, we get 
\begin{equation}
  \label{eq:rho-t-sl2}
  \frac{d\rho}{dt} (te) = \sqrt{(2k)^4+\frac{4c}{t^2}}.
\end{equation}

\begin{proposition}
  \label{prop:unique}
  The Kähler potential~$\rho$ of Proposition~\ref{prop:sl2-potentials} is a
  hyperKähler potential if and only if $c=0$.  In this case,
  $\rho=2k\sqrt\eta$ and $\rho(te)=4k^2t$. 
\end{proposition}

\begin{proof}
  Let $Y=(d\rho)^\sharp$ be the vector field dual to~$d\rho$.  If $\rho$ is
  a hyperKähler potential then $IY$ is an isometry preserving~$I$ and
  $\rho$ is the corresponding moment
  map~\cite[Proposition~5.5]{Swann:MathAnn}.  Now $d\rho=\rho'd\eta
  =2k^2\rho'\re\inp\cdot{\sigma X}$, whereas, by Remark~\ref{rem:metric},
  \begin{equation*}
    g(Y,\xi_A) = 2\re(\rho'k^2\inp{\xi_A}{\sigma
    Y}+\rho''k^4\inp{\xi_A}{\sigma X}\inp X{\sigma Y}).
  \end{equation*}
  So we have
  \begin{equation*}
    \rho' X = \rho' Y +\rho''k^2\inp Y{\sigma X} X,
  \end{equation*}
  which implies that $Y=\lambda X$ with 
  \begin{equation*}
    \begin{split}
      \lambda &= \frac{\rho'}{\rho'+\eta\rho''}
      = \frac{2\eta {\rho'}^2}{k^2} = 2 + \frac{2c}{k^2\eta},
    \end{split}
  \end{equation*}
  using~\eqref{eq:ode}.  The vector field $X$ is generated by scaling in
  the nilpotent orbit, so $IX$ preserves the complex structure~$I$.  Now
  \begin{equation*}
    \begin{split}
      (L_{IY}I)(Z) &= [\lambda IX,IZ] - I [\lambda IX, Z] \\
      &= \lambda(L_{IX}I) Z - ((IZ)\lambda) IX + (Z\lambda) X,
  \end{split}
  \end{equation*}
  so with $L_{IX}I=0$, we have $L_{IY}I=0$ only if $\lambda$~is constant.
  But this is exactly the requirement that $c=0$.
\end{proof}

\begin{remark}
  \label{rem:Eguchi-Hanson}
  The substitution $k^2\eta+c=(\frac r2)^4$ in equation~\ref{eq:metric}
  shows that these are the Eguchi-Hanson metrics
  (cf.~\cite{Dancer:survey-hK}).  See~\cite{Kobak-Swann:Wolf} for details.
\end{remark}

Let us now turn to the regular nilpotent orbit~$\mathcal O_\Delta$
in~$\so(4,\mathbb C)$.  As in the proof of Theorem~\ref{thm:so4}, we write
$\so(4,\mathbb C)=\Sl(2,\mathbb C)_+ \oplus \Sl(2,\mathbb C)_-$ and note
that $\mathcal O_\Delta=\mathcal O_+\times\mathcal O_-$ where $\mathcal
O_\pm$ is the nilpotent orbit in~$\Sl(2,\mathbb C)_\pm$.  Let
$\{e_\pm,f_\pm,h_\pm\}$ be bases for~$\Sl(2,\mathbb C)_\pm$ as
in~\eqref{eq:efh}.  Again we will use the inner product which is
$k^2\inp[\so(4)]\cdot\cdot$.

Using the action of~$\SO(4)$, we may take our representative element~$X$
of~$\mathcal O_\Delta$ to be $X=X_++X_-= s\, e_+ + t\, e_-$ with $s,t >0$.
We have one invariant for each $\Sl(2,\mathbb C)$: we write $\eta_\pm =
k^2\inp[\Sl(2)]{X_\pm}{\sigma X_\pm}$, so $\eta_+=4k^2s^2$, etc.  Let
$\rho_+=\partial\rho/\partial\eta_+$, etc.  Then we may calculate the
Kähler form $\omega_I$ and the candidate almost complex structure $J$ as
in~\S\ref{sec:2-i}.  For the Kähler form we get
\begin{equation*}
  \begin{split}
    \omega_I(\xi_A,\xi_B)
    = 2k^2\im\Bigl(&
    \rho_+ \inp{\xi_A^+}{\sigma\xi_B^+}
    +
    \rho_- \inp{\xi_A^-}{\sigma\xi_B^-}
    \\
    &\quad +
    \rho_{++} k^2 \inp{\xi_A^+}{\sigma X_+}
    \inp{\sigma\xi_B^+}{X_+}\\ 
    &\quad +
    \rho_{+-} k^2
    \bigl(
    \inp{\xi_A^+}{\sigma X_+} \inp{\sigma\xi_B^-}{X_-}\\
    &\quad\hphantom{+ \rho_{+-} k^2}\quad +
    \inp{\xi_A^-}{\sigma X_-} \inp{\sigma\xi_B^+}{X_+}
    \bigr) \\
    &\quad +
    \rho_{--} k^2 \inp{\xi_A^-}{\sigma X_-}
    \inp{\sigma\xi_B^-}{X_-} 
    \Bigr),
  \end{split}
\end{equation*}
where $\xi_A^+=[A,X_+]$, etc.  The endomorphism $J$ is given by
\begin{multline*}
  J_X(\xi_A^+)
  = - 2\rho_+ [X_+,\sigma\xi_A^+]\\
  - 2k^2 \inp{\sigma\xi_A^+}{X_+} \bigl( \rho_{++} [X_+,\sigma X_+] +
  \rho_{+-} [X_-,\sigma X_-] \bigr)
\end{multline*}
and a similar expression for $\xi_A^-$.  In particular, $J_X h_+ = - 4 s\,
\rho_+ e_+$ and
\begin{equation*}
  J_X e_+ = 2 s (\rho_+ + \eta_+\rho_{++}) h_+ + 2 \frac{t^2}s \eta_+
  \rho_{+-} h_-.
\end{equation*}
Thus the $\Sl(2,\mathbb C)_-$-component of $J^2_X h_+$ is a constant times 
$\eta_+\eta_-\*\rho_+\rho_{+-}\*h_-$.  For $J_X^2$ to be~$-1$, we need
$\rho_+\rho_{+-}=0$, which implies $\partial(\rho_+^2)/\partial \rho_- =0$
and hence $\rho_{+-}=0$.  Thus $J_X$ preserves the $\Sl(2,\mathbb
C)$-summands of~$\so(4,\mathbb C)$.

\begin{proposition}
  \label{prop:so4-potentials}
  Any hyperKähler structure on the regular orbit \linebreak $\mathcal
  O_\Delta=\mathcal O_+\times\mathcal O_-$ of~$\so(4,\mathbb C)$ which is
  $\SO(4)$-invariant, admits a Kähler potential and has complex-symplectic
  form $k^2\omega_c^{\mathcal O_\Delta}$, is a product of
  $\SU(2)$-invariant structures on the factors~$\mathcal O_\pm$, and these
  are given by Proposition~\ref{prop:sl2-potentials}. \qed
\end{proposition}

\section{Potentials for Next-to-Minimal Orbits}
\label{sec:next-to-min}
We now come to the main result of this paper.  We consider next-to-minimal
orbits with compatible $G$-invariant hyperKähler metrics, except for
$G=\SU(3)$.  We show that such metrics admitting a hyperKähler potential are
unique, and we calculate the potential.

If we assume that the potential is only Kähler, we still have uniqueness in
some cases, but we get a list of exceptions: orbits which admit a
one-parameter family of hyperKähler metrics. These can be thought of as a
generalisation of the Eguchi-Hanson metric (cf.\ 
Remark~\ref{rem:Eguchi-Hanson}).

\begin{theorem}
  \label{thm:main}
  Suppose $G$ is a compact simple Lie group and $\mathcal O$ is a nilpotent
  orbit in~$\mathfrak g^{\mathbb C}$ of cohomogeneity two.  
  \begin{enumerate}
  \item [(i)] $\mathcal O$ admits a unique $G$-invariant compatible
    hyperKähler metric with hyperKähler potential. This potential is given
    by
    \begin{equation}
      \label{eq:rho-generic}
      \rho = 2k\sqrt{\eta_1 + 2\sqrt{\tfrac12\eta_1^2 - k^2\eta_2}}
    \end{equation}
    for $\mathfrak g\ne\lieG_2$, where the constant $k$ is given in
    Table~\ref{tab:k2}, and, for $\lieG_2$,
    \begin{equation}
      \label{eq:rho-g2}
      \rho = \sqrt8\sqrt{\eta_1 + \sqrt6\sqrt{\eta_1^2 - 4\eta_2}}.
    \end{equation}
  \item [(ii)] The above metric on $\mathcal O$ is in fact a unique
    $G$-invariant compatible hyperK\"ahler metric with a \emph{Kähler}
    potential unless $\mathfrak g$ is one of $\LSP(2)\cong\so(5)$,
    $\su(4)\cong\so(6)$, $\so(8)$ or $\mathcal O$~is of Jordan
    type~$(31^{n-3})$ in~$\so(n)$.  In these cases, the metric lies in a
    one-parameter family of hyperK\"ahler metrics with Kähler potentials.
  \end{enumerate}
\end{theorem}

\begin{table}[tbp]
  \begin{center}
    \newcommand{\tablestrut}{\vrule height 14pt depth 6pt width 0pt}
    \leavevmode
    \begin{tabular}[t]{c||c|c|c|c|c|c}
      \tablestrut Type&$A_n$, $C_n$&$B_n$,
      $D_n$&$\LieF_4$&$\LieE_6$&$\LieE_7$&$\LieE_8$\\
      \hline
      \tablestrut $k^2$&$\tfrac 12(n+1)$&$\tfrac 12(n-2)$&$\tfrac
      92$&$6$&$9$&$\tfrac{35}2$\\
    \end{tabular}
    \caption{The constant $k^2$ in the potentials of Theorem~\ref{thm:main}.}
    \label{tab:k2}
  \end{center}
\end{table}

\begin{remark}
  Note that the Theorem provides hyperKähler potentials for all
  next-to-minimal orbits, except when $\mathfrak g=\su(3)$.  However, the
  potential in this remaining case was computed
  in~\cite{Kobak-Swann:nilpotent}, see
  also~\cite{Kobak-Swann:finite-potentials}.
\end{remark}

We divide the proof of the Theorem into three parts.

\subsection{The General Case}
\label{sec:general}
This is when $\mathfrak g$ is neither $\su(3)$ nor~$\lieG_2$.

Let $X$~be a generic element of~$\mathcal O$.  By Theorem~\ref{thm:so4},
$X$~lies in the regular orbit~$\mathcal O_\Delta$ of a real $\so(4,\mathbb
C)$-subalgebra.

For $\rho(\eta_1,\eta_2)$ to be a hyperKähler potential for~$\mathcal O$ it
is necessary that $\rho$~is a Kähler potential for an invariant hyperKähler
structure on~$\mathcal O_\Delta$.  To see this, first note that
equation~\eqref{eq:dId} is invariant by pull-back under the inclusion map
$\mathcal O_\Delta \hookrightarrow \mathcal O$.  Now equation~\eqref{eq:J}
shows that $J\xi_A$ remains in the subalgebra generated by $A$, $X$ and
$\sigma X$.  Thus if $A\in\so(4,\mathbb C)$, so is $JA$ and thus $\mathcal
O_\Delta$ is a hyperKähler submanifold of~$\mathcal O$.

As in~\S\ref{sec:models}, write $\so(4,\mathbb C)=\Sl(2,\mathbb
C)_+\oplus\Sl(2,\mathbb C)_-$, $\mathcal O_\Delta=\mathcal
O_+\times\mathcal O_-$ and $X=X_++X_-=se_++te_-$.  Our two invariants
on~$\mathcal O$ are given by
\begin{equation*}
  \begin{split}
    \eta_1(X)
    & = \inp[\mathfrak g] X{\sigma X} = - \inp[\mathfrak g]{se_+ +
    te_-}{sf_+ + tf_-}\\ 
    & = -(s^2+t^2) k^2 \inp[\su(2)] ef = 4k^2 (s^2+t^2),
  \end{split}
\end{equation*}
and a similar computation gives
\begin{equation*}
  \eta_2(X) = 8k^2 (s^4+t^4),
\end{equation*}
where $k^2$~is the constant such that $\inp[\mathfrak g]\cdot\cdot
|_{\so(4,\mathbb C)} = k^2 \inp[\so(4)]\cdot\cdot$.  Now
\begin{equation*}
  \begin{split}
    d\rho
    &= \rho_1d\eta_1+\rho_2d\eta_2 \\
    &= 8k^2
    \left(
      s(\rho_1+4s^2\rho_2)ds
      +
      t(\rho_1+4t^2\rho_2)dt
    \right),
  \end{split}
\end{equation*}
so $\rho_s := \partial \rho/\partial s = 8k^2s(\rho_1+4s^2\rho_2)$,
etc., and solving for $\rho_1$ and $\rho_2$ we get
\begin{equation*}
  \rho_1 = -\frac{t^3\rho_s - s^3\rho_t}{8k^2st(s^2-t^2)},
  \qquad
  \rho_2 = \frac{t\rho_s-s\rho_t}{32k^2st(s^2-t^2)}.
\end{equation*}
Note that, by Proposition~\ref{prop:so4-potentials}
and~\eqref{eq:rho-t-sl2}, ${\rho_s}^2=16k^4+c_+/s^2$ and
${\rho_t}^2=16k^4+c_-/t^2$, for some constants~$c_\pm$.

The elements $X_+$ and $X_-$ lie in the closure of~$\mathcal O_\Delta$ and
hence of~$\mathcal O$; so $X_\pm$ lie in the minimal nilpotent orbit
of~$\mathfrak g^{\mathbb C}$.  We deduce that $M_+:= G/N(\SU(2)_+)$ is a
Wolf space and hence, since $\SU(2)_+$~corresponds to a highest
root~\cite{Wolf:quaternionic},
\begin{equation}
  \label{eq:gC-dp}
  \mathfrak g^{\mathbb C} = \Sl(2,\mathbb C)_+ + \mathfrak k_+ + S^1_+\otimes
  E_+,
\end{equation}
where $\mathfrak k_+$ commutes with $\Sl(2,\mathbb C)$, $E_+$~is a
non-trivial representation of~$\mathfrak k_+$ and $S^1_+\cong\mathbb
C^2$~is the fundamental representation of $\Sl(2,\mathbb C)_+$.  On the
other hand, we have a similar decomposition of~$\mathfrak g^{\mathbb C}$
corresponding to~$\Sl(2,\mathbb C)_-$.  As $\Sl(2,\mathbb C)_+$ and
$\Sl(2,\mathbb C)_-$ commute with each other, we deduce that $\Sl(2,\mathbb
C)_-\subset \mathfrak k_+$ and that $E_+\supset S^1_-$.  So as an
$\so(4,\mathbb C)$-module, $\mathfrak g^{\mathbb C}$ always contains a copy
of~$S^1_+\otimes S^1_-$.

On the orthogonal complement to~$\so(4,\mathbb C)$, we have,
from~\eqref{eq:J},
\begin{equation}
  \label{eq:J-so4-p}
  \begin{split}
    J_X \xi_A &= - 2\rho_1[X,\sigma \xi_A]\\
    &\quad+ 4\rho_2
    \left(
      2[X,[\sigma X,[X,\sigma \xi_A]]]
      - [X,[X,[\sigma X,\sigma \xi_A]]]
    \right).
\end{split}
\end{equation}

First suppose that $E_+$ contains a trivial $\Sl(2,\mathbb
C)_-$-module~$\mathbb C^r$; take $r$~maximal.  The real structure~$\sigma$
preserves the module $S^1_+\otimes\mathbb C^r$ and acts on
$S^1_+\cong\mathbb H$ as~$j$, so $\mathbb C^r$ has a quaternionic
structure~$\mathfrak j$ and is even-dimensional.  Choose a basis for
$S^1_+$ so that $\ad e_+$ acts as $\begin{spmatrix} 0&1\\0&0
\end{spmatrix}$.  Then any tangent vector $\xi_A\in
S^1_+\otimes\mathbb C^r$ has the form $\begin{spmatrix}
1\\0\end{spmatrix}\otimes v$ and we have
\begin{equation*}
  J_X\xi_A = -2s(\rho_1+4s^2\rho_2)
    \begin{spmatrix}
      1\\0
    \end{spmatrix}
  \otimes \mathfrak j v
  = - \tfrac1{4k^2}\rho_s 
    \begin{spmatrix}
      1\\0
    \end{spmatrix}
  \otimes \mathfrak j v
  .
\end{equation*}
Thus $J^2=-1$ on~$S^1_+\otimes\mathbb C^r$ if and only if
${\rho_s}^2=16k^4$.  This implies that the constant $c_+$ is zero if $E_+$
has an trivial $\Sl(2,\mathbb C)_-$-submodule.

The existence of an trivial $\Sl(2,\mathbb C)_-$-submodule in $E_+$ is not
guaranteed.  However, we do always have an $S^1_-$-summand, so we now
consider the case when $\xi_A$ lies in an $\so(4)$-module $S^1_+\otimes
S^1_-$.  This is Killing orthogonal to~$\so(4,\mathbb C)$.  We choose bases
so that $\ad X$ acts as
\begin{equation*}
  s 
    \begin{spmatrix}
      0&1\\
      0&0
    \end{spmatrix}
  \otimes\Id
  + t \Id \otimes
  \begin{spmatrix}
    0&1\\
    0&0
  \end{spmatrix}
\end{equation*}
and $\sigma=j\otimes j$ for $j$~the standard quaternionic structure
on~$S^1\cong\mathbb H$.  The image of $\ad X$ is two-dimensional and
spanned by
\begin{equation*}
  \xi_1 := 
  \begin{spmatrix}
    1\\0
  \end{spmatrix}
  \otimes
  \begin{spmatrix}
    1\\0
  \end{spmatrix}
  \quad\text{and}\quad
  \xi_2 := s 
  \begin{spmatrix}
    1\\0
  \end{spmatrix}
  \otimes
  \begin{spmatrix}
    0\\1
  \end{spmatrix}
  + t
  \begin{spmatrix}
    0\\1
  \end{spmatrix}
  \otimes
  \begin{spmatrix}
    1\\0
  \end{spmatrix}
  .
\end{equation*}
These satisfy 
\begin{gather*}
  [X,\xi_1]=0, \qquad [\sigma X,\xi_1]=\sigma \xi_2,\\
  [X,\xi_2]=2st\xi_1\quad\text{and}\quad
  [\sigma X,\xi_2]= -(s^2+t^2)\sigma \xi_1.
\end{gather*}
So, equation~\eqref{eq:J-so4-p} gives
\begin{gather*}
  J\xi_1 = -2 (\rho_1 +4\rho_2 (s^2+t^2)) \xi_2,\\
  J\xi_1 = 2(\rho_1 (s^2+t^2) + 4\rho_2 (s^4+t^4)) \xi_1.
\end{gather*}
Substituting for $\rho_1$ and $\rho_2$ in terms of $\rho_s$ and $\rho_t$,
gives 
\begin{equation*}
  J^2\xi_1 = - \frac{t^2{\rho_t}^2-s^2{\rho_s}^2}{16k^4(t^2-s^2)} \xi_1.
\end{equation*}
So $J^2=-1$ on $S^1_+\otimes S^1_-$ if and only if
$t^2{\rho_t}^2-s^2{\rho_s}^2 = 16k^4(t^2-s^2)$.  But ${\rho_s}^2=
16k^4+c_+/s^2$, etc., so $c_+=c_-$.

We conclude that if $E_+$ contains a trivial $\Sl(2,\mathbb C)_-$-summand,
then $c_+=c_-=0$.  This gives $\rho_s = 4k^2$ and $\rho_t=4k^2$, so
$\rho(s,t)=4k^2(s+t)$.  Rewriting this in terms of $\eta_1$ and $\eta_2$
gives the potential in the Theorem.  If $E_+$ does not have a trivial
summand, we get a one-parameter family of potentials and hyperKähler
metrics with $c_+=c_-$.

It remains to determine the constant~$k$ and when $E_+$~contains a trivial
$\SU(2)_-$-module.  The decomposition~\eqref{eq:gC-dp} gives the action
of~$\ad e_+$ and hence the Killing inner product $\inp[\mathfrak
g]{e_+}{\sigma e_+}$ is $4+\dim_{\mathbb C}E_+$, since
$\inp[\su(2)_+]{e_+}{\sigma e_+}=4$.  So $k^2=(4+\dim_{\mathbb C}E_+)/4$.
Moreover, $S^1_+\otimes E_+=TM_+\otimes\mathbb C$, so $\dim_{\mathbb
C}E_+$~is half the real dimension of the Wolf space~$M_+$, which may be
found in, e.g., Besse~\cite[p.~409]{Besse:Einstein}, or read-off from the
discussion below.  This leads to Table~\ref{tab:k2}.

Finally, we determine the decompositions of~$E_+$ under the action
of~$\Sl(2,\mathbb C)_-$.

If $G=\SU(n)$, then $\mathfrak k_+\cong\un(n-2)$, and $E_+=\mathbb C^{n-2}$
is the fundamental representation twisted by a representation of the
central~$\un(1)$.  Now $\Sl(2,\mathbb C)_-$ corresponds to a highest root
vector in~$\mathfrak k_+$, so $E_+=S^1_-+\mathbb C^{n-4}$ as a
$\Sl(2,\mathbb C)_-$-module.  So for $n=4$, we have a one-parameter family
of potentials $c_+=c_-$, and for $n>4$, the potential is unique.

For $G=\SP(n)$, $\mathfrak k_+\cong\LSP(n-1,\mathbb C)$ and $E_+\cong
\mathbb C^{2n-2}\cong \mathbb H^{n-1}$~is the fundamental representation.
Under the highest root $\Sl(2,\mathbb C)$, this representation splits as
$S^1_-+\mathbb C^{2n-4}$, so for $n>1$, we have a unique potential.

In the case $G=\SO(n)$, there are two orbit types to consider.  The
centraliser $\mathfrak k_+= \Sl(2,\mathbb C)+\so(n-4,\mathbb C)$ and there
are two choices for $\Sl(2,\mathbb C)_-$, one in each summand of~$\mathfrak
k_+$.  When $\Sl(2,\mathbb C)_-=\Sl(2,\mathbb C)$, we get $E_+\cong
S^1_-\otimes \mathbb R^{n-4}$, and there is a one-parameter family of
potentials.  On the other hand, if $\Sl(2,\mathbb C)_-$ lies in the
summand~$\so(n-4,\mathbb C)$, then $E_+\cong \mathbb
C^2\otimes(S^1_-+\mathbb R^{n-8})$.  For $n>8$, this gives a unique
potential, but for $n=8$, we again get a family.

We now come to the four exceptional cases.  Firstly, if $G=\LieF_4$, then
$\mathfrak k_+\cong\LSP(3,\mathbb C)$ and if $E=\mathbb H^3$ is the
fundamental representation, then $E_+ \cong \Lambda_0^3E = \Lambda^3 E -
E$, is a $14$-dimensional irreducible representation.  For a highest root
$\Sl(2,\mathbb C)_-$ in~$\LSP(3,\mathbb C)$, we have $E\cong S^1_-+\mathbb
C^4$ and hence $E_+ \cong 5S^1_-+\mathbb C^4$.  So $E_+$ has a trivial
summand and hence the potential is unique.

For $G=\LieE_6$, $\mathfrak k_+ \cong \Sl(6,\mathbb C)$ and $E_+\cong
\Lambda^{3,0}\mathbb C^6$.  Under a highest root $\Sl(2,\mathbb C)_-$,
we have $\Lambda^{1,0}\mathbb C^6\cong S^1_-+\mathbb C^4$ and hence $E_+ =
4S^1_-+\mathbb C^8$, giving a unique potential.

When $G=\LieE_7$, $\mathfrak k_+ \cong \so(12,\mathbb C)$ and $E_+\cong
\Delta^{12}_+$, the positive spin representation.  For a highest root
$\Sl(2,\mathbb C)_-$, the normaliser in~$\so(12,\mathbb C)$ is
$\Sl(2,\mathbb C)_-+\Sl(2,\mathbb C)+\so(8,\mathbb C)$ and the fundamental
representation of~$\SO(12)$ decomposes as $\mathbb C^{12}\cong S^1_-\otimes
S^1 + V$, where $V\cong\mathbb C^8$ is the fundamental representation
of~$\so(8,\mathbb C)$.  The spin representation splits as
$\Delta^{12}_+\cong S^1_-\otimes\Delta^8_+ + S^1\otimes\Delta^8_-$, and so
$E_+\cong 8S^1_-+\mathbb C^{16}$ has a trivial summand.

Finally, for $G=\LieE_8$, $\mathfrak k_+\cong \lieE_7^{\mathbb C}$ and
$E_+\cong {\scriptstyle 100\Dyatop0000}$.  A highest root
$\Sl(2,\mathbb C)_-$ in~$\lieE_7^{\mathbb C}$ has
centraliser~$\so(12,\mathbb C)$ and $E_+\cong 12 S^1_- + \mathbb C^{32}$,
where $\mathbb C^{32}\cong\Delta^{12}_+$.  So again we get a unique potential.

\subsection{The Exceptional Case~$\LieG_2$}
\label{sec:G2}
The Dynkin diagram for the next-to-minimal orbit~$\mathcal O$ in~$\LieG_2$
is $\scriptstyle 0\rthreea1$.  This says that there is a basis
$\{\alpha,\beta\}$ for the simple positive roots, with $\alpha$~short and
$\beta$~long, such that $\ad h$ acts on $\mathfrak g_\alpha$ and $\mathfrak
g_\beta$ with eigenvalues $1$ and~$0$ respectively.  We thus have
$\mathfrak g(2)=\mathfrak g_{\beta+2\alpha}$ and $\mathfrak g(3)=\mathfrak
g_{\beta+3\alpha}\oplus\mathfrak g_{2\beta+3\alpha}$.  From the discussion
in~\S\ref{sec:subalgebras}, the isotropy group~$\SU(2)\Un(1)$ of the
Beauville bundle acts transitively on the unit sphere in~$\mathfrak g(3)$,
so using the action of the compact group~$\LieG_2$, we can move a typical
element of~$\mathcal O$ to $X\in\mathfrak g_{\beta+2\alpha}\oplus\mathfrak
g_{2\beta+3\alpha}$.  We may thus write $X=s E_{\beta+2\alpha}+t
E_{2\beta+3\alpha}$, with $s,t>0$, where $E_i$ are such that for
$F_i:=-\sigma E_i$ and $H_i=[E_i,F_i]$ we have $[H_i,E_i]=2E_i$.

At $X$, our two invariants are
\begin{equation*}
  \eta_1(X) = 8(s^2+3t^2)
  \quad\text{and}\quad
  \eta_2(X) = 16(s^4+6s^2t^2+3t^4).
\end{equation*}
As in the previous section, we compute $J^2$ on particular tangent vectors
using~\eqref{eq:J} and then rewrite the equations in terms of $s$ and~$t$.
This is quite hard work to do by hand, and so we used \textsc{Maple} to do
the following computations.  The code for this is described
in~\cite{Kobak-Swann:G2-Maple}.  

On $\mathfrak g_{\alpha+\beta}$, one finds that $J^2=-1$ only if
\begin{equation}
  \label{eq:rho-g2-1}
  \frac1{64s}\rho_s(s\rho_s+t\rho_t) = 1,
\end{equation}
where $\rho_s$ is $\partial \rho/\partial s$, etc.  Now
$X=[X,H_\beta-H_\alpha]=[X,3H_{2\beta+3\alpha}-5H_{\beta+2\alpha}]$, so
$X$~is tangent to the orbit~$\mathcal O$.  The condition $J^2X=-X$, gives
the following three equations
\begin{subequations}
  \begin{gather}
    \label{eq:J2X-a}
    \begin{split}
      s (2 s \rho_s + t \rho_t) \rho_{ss} + t (t \rho_t + 3 s \rho_s)
      \rho_{st} + t^2 \rho_s \rho_{tt} \qquad&\\
      + 2 (t \rho_t + s \rho_s)\rho_s &= 128 s,
    \end{split}
    \\
    \label{eq:J2X-b}
    9 s \rho_s \rho_{ss} + (9 t \rho_s + s \rho_t)
    \rho_{st} + t \rho_t \rho_{tt} + 9 \rho_s^2 + \rho_t^2=576\\
    \label{eq:J2X-c}
    \begin{split}
      3 s t (9 t \rho_s + s \rho_t) \rho_{ss} - s t (s \rho_t - 3 t
      \rho_s) \rho_{tt}\quad&\\
      + (3t(s^2 + 9t^2) \rho_s + s(3t^2 - s^2) \rho_t) \rho_{st} &=
      (s \rho_t - 9 t \rho_s) (s \rho_t + 3 t \rho_s)
    \end{split}
  \end{gather}
\end{subequations}
by considering the components in $\mathfrak g_{\beta+2\alpha}$, $\mathfrak
g_{2\beta+3\alpha}$ and $\mathfrak g_{2\beta+\alpha}$.  Considering
$9s\rho_s$ times~\eqref{eq:J2X-a} minus $s (2 \rho_s s + t \rho_t)$
times~$~\eqref{eq:J2X-b}$ gives a new equation not involving $\rho_{ss}$.
In a similar way, we may eliminate~$\rho_{ss}$ form the pair of equations
\eqref{eq:J2X-a} and~\eqref{eq:J2X-c}.  Eliminating $\rho_{st}$ from these
two new equations not involving~$\rho_{ss}$, we get the following equation
which does not involve~$\rho_{tt}$:
\begin{equation*}
  s^3t(2s\rho_s+t\rho_t)(s\rho_t-9t\rho_s)^2=0.
\end{equation*}
Thus either 
\begin{equation}
  \label{eq:rho-g2-2}
  \text{(i)}\quad\rho_t = - 2\frac st \rho_s,
  \qquad\text{or}\qquad
  \text{(ii)}\quad\rho_t = 9\frac ts \rho_s.
\end{equation}
In case~(i), substituting into \eqref{eq:rho-g2-1} one gets $\rho_s^2=-64$,
which has no (real) solutions.  In case~(ii), we have 
\begin{equation*}
  \rho_s = \varepsilon \frac{8s}{\sqrt{s^2+9t^2}},\qquad
  \rho_t = \varepsilon \frac{72t}{\sqrt{s^2+9t^2}},
\end{equation*}
where $\varepsilon\in\{\pm1\}$.  Integrating we find that
\begin{equation}
  \label{eq:rho-g2-f}
  \rho = \varepsilon 8\sqrt{s^2+9t^2}.
\end{equation}
To get a positive-definite metric, take~$\varepsilon=+1$.  Rewriting
\eqref{eq:rho-g2-f} in terms of~$\eta_1$ and~$\eta_2$ gives the claimed
result.  One may check directly that the resulting $J$ satisfies~$J^2=-1$
on the whole tangent space.

\subsection{Uniqueness of HyperKähler Potentials}
\label{sec:hk-potentials}
The only statement left to verify in the proof of Theorem~\ref{thm:main},
is that equations \eqref{eq:rho-generic} and~\eqref{eq:rho-g2} give the
unique compatible hyperKähler potentials on the orbits.  In the cases, when
the Kähler potential is unique there is nothing to prove, because the
general theory~\cite{Swann:MathAnn} gives the existence of such a
potential.  We may therefore assume we are in the general case and that our
generic element~$X$ lies in a real $\so(4,\mathbb C)$-subalgebra.  Now
$\mathcal O_\Delta$ is a hyperKähler submanifold of~$\mathcal O$, and so
by~\eqref{eq:dId}, $\rho$ is a hyperKähler potential for~$\mathcal O$ only
if it restricts to a hyperKähler potential for~$\mathcal O_\Delta$.
However, the hyperKähler structure on~$\mathcal O_\Delta$ is the product of
two hyperKähler structures on~$\Sl(2,\mathbb C)$-orbits and on each of
these factors the hyperKähler potential is unique by
Proposition~\ref{prop:unique}.  Thus there is only one hyperKähler
potential compatible with the structure of~$\mathcal O$. \qed

\begin{remark}
  The hyperKähler metrics constructed in Theorem~\ref{thm:main} have an
  extra $\Un(1)$-symmetry given by $X\mapsto e^{i\theta}X$ which preserves
  the complex structure~$I$ but moves~$J$.  In the case of~$\Sl(2,\mathbb
  C)$, the metrics are of Bianchi type~$IX$ and it is known, e.g.,
  from~\cite{Belinskii-GPP:Bianchi-IX}, that there are triaxial hyperKähler
  metrics that do not have $\Un(2)$-symmetry.  Thus concentrating on
  metrics admitting a Kähler potential is a genuine restriction.
\end{remark}

\begin{remark}
  The one-parameter families in Theorem~\ref{thm:main} occur exactly when
  $E_+\cong \mathbb C^s\otimes S^1_-$.  Considering the weights of the
  action of a semi-simple element in the diagonal $\Sl(2,\mathbb
  C)$-subalgebra of~$\so(4,\mathbb C)$ on~$\mathfrak g^{\mathbb C}$, we see
  that this exactly the case when $\mathfrak g(1)=0$.  This says that the
  Beauville bundle coïncides with the cotangent bundle~$T^*\mathcal F$,
  rather than being a proper subbundle.  In the case of the one-parameter
  families $\mathcal F=\Gro_2(\mathbb R^n)$ and for $c\ne0$, the Kähler
  potentials extend to give non-singular metrics on~$T^*\mathcal F$,
  generalising the Eguchi-Hanson metrics on~$T^*\CP(1)$.  As $\mathcal
  F$~is Hermitian symmetric, this is one of the cases considered by Biquard
  \&\ Gauduchon \cite{Biquard-Gauduchon:potential}.
\end{remark}

\providecommand{\bysame}{\leavevmode\hbox to3em{\hrulefill}\thinspace}

\end{document}